\documentclass[11pt]{article}

\usepackage{amsfonts}
\usepackage{amssymb}
\usepackage{amsmath}
\usepackage{amsthm}

\def\negthickspace{\!\!\!}

\newcommand{\nicefrac}[2]
{\leavevmode \kern.1em\raise.5ex\hbox{\the\scriptfont0 #1}
             \kern-.1em/\kern-.15em\lower.25ex
             \hbox{\the\scriptfont0 #2}}

\newtheorem*{theorem}{Theorem}
\newtheorem*{proposition}{Proposition}
\newtheorem*{corollary}{Corollary}

\newtheorem*{definition}{Definition}
\newtheorem*{lemma}{Lemma}

\theoremstyle{definition}
\newtheorem*{remark}{Remark}
\newtheorem*{remarks}{Remarks}

\theoremstyle{definition}
\newtheorem*{example}{Example}
\newtheorem*{examples}{Examples}

\begin{document} 

\begin{center}
{\Large{\sc Area estimates}}\\[0.2cm]
{\Large{\sc for two-dimensional immersions of mean curvature type}}\\[0.2cm]
{\Large{\sc in Euclidean spaces of higher codimension}}\\[1cm]
{\large Steffen Fr\"ohlich}\\[0.4cm]
{\small\bf Abstract}\\[0.4cm]
\begin{minipage}[c][2.5cm][l]{12cm}
{\small We establish area bounds for two-dimensional immersions in $\mathbb R^3$ and $\mathbb R^n.$ Namely, for $\mu$-stable immersions in $\mathbb R^3$ ($\mathbb R^n$), for graphs in $\mathbb R^3$ which solve quasilinear equations in divergence form, and for graphs which are critical for Fermat-type variational problems in $\mathbb R^n.$}
\end{minipage}
\end{center}
{\small MCS 2000: 35J60, 53A07, 53A10}
\section{Introduction}
\subsection{The results}
In this paper we prove area bounds for the following types of surfaces: 
\begin{itemize}
\item[1.]
$\mu$-stable immersions of prescribed mean curvature-type in $\mathbb R^3$ ($\mathbb R^n$) in terms of a suitable stability constant and of curvature terms (chapter 2, section 2.3, 2.5, 2.6);
\vspace*{-0.6ex}
\item[2.]
graphs $(x,y,\zeta(x,y))$ of mean curvature type in $\mathbb R^3,$ which are solutions of non-homogeneous divergence form equations (chapter 3, section 3.2);
\vspace*{-0.6ex}
\item[3.]
graphs $(x,y,\zeta_1(x,y),\ldots,\zeta_{n-2}(x,y))$ in $\mathbb R^n,$ $n\ge 3,$ which are critical for Fermat-type variational problems (chapter 4, section 4.3).
\end{itemize}
For example, area bounds are crucial for compactness results (see e.g. \cite{Clarenz_Mosel_01} for such results concerning weighted minimal surfaces, see section 2.2 below), but also for various gradient and curvature estimates for nonlinear differential systems: Here, we mention to \cite{Bergner_Froehlich_01}, where such a bound for Fermat-type graphs in $\mathbb R^n$ was left unproved, further \cite{Froehlich_01}-\cite{Froehlich_05}, \cite{Sauvigny_01}, \cite{Sauvigny_02}, where area bounds are used for various curvature estimates, and, finally, \cite{Froehlich_Winklmann_01} for curvature estimates of $n$-manifolds in $\mathbb R^m$ using techniques from \cite{Ecker_01}, \cite{Schoen_Simon_Yau_01}.\\[1ex]
We will not discuss isoperimetric inequalities. But we mention that such inequalities are established e.g. in \cite{Clarenz_Mosel_01} for immersed critical points of elliptic variational problems (\ref{2.8}) using Fourier series methods, or in \cite{Ecker_01} (and references therein) for mean curvature immersions using a generalized Sobolev inequality. We want to extend this later method in a further paper to prove other area bounds for immersions of mean-curvature type (see the discussion in section 2.2).
\section{$\mu$-stable geodesic discs of mean curvature-type}
\setcounter{equation}{0}
In the first part of this note we consider immersions
\begin{equation}\label{2.1}
  X=X(u,v)=\big(x^1(u,v),x^2(u,v),x^3(u,v)\big)\in C^3(B,\mathbb R^3)
\end{equation}
on the closed unit disc $B:=\big\{(u,v)\in\mathbb R^2\,:\,u^2+v^2\le 1\big\}$ such that $\mbox{rank}\,\partial X=2$ in $B$ for its Jacobian $\partial X\in\mathbb R^{3\times 2}.$ The unit normal vector of $X$ is defined as
\begin{equation}\label{2.2}
  N=\frac{X_u\times X_v}{|X_u\times X_v|}\quad\mbox{in}\ B
\end{equation}
with the partial derivatives $X_u$ and $X_v$ of $X,$ and $\times$ means the usual vector product in $\mathbb R^3.$\\[1ex]
In the next two sections we specify the class of immersion we will deal with.
\subsection{First step: Introduction of weighted metrics}
We equip the immersions with a weighted metric of Finsler type: Let us given a symmetric and positive definite weight matrix
\begin{equation}\label{2.3}
  {\mathbf G}(X,Z)\in C^2(\mathbb R^3\times\mathbb R^3\setminus\{0\},\mathbb R^{3\times 3})
\end{equation}
with the properties: For all $(X,Z)\in\mathbb R^3\times\mathbb R^3\setminus\{0\}$ there hold
\begin{itemize}
\item[(G1)]
${\mathbf G}(X,Z)={\mathbf G}(X,\lambda Z)$ for all real $\lambda>0;$
\vspace*{-0.6ex}
\item[(G2)]
${\mathbf G}(X,Z)\circ Z^t=Z^t,$ where the upper $t$ denotes transposition;
\vspace*{-0.6ex}
\item[(G3)]
$(1+g_0)^{-1}|\xi|^2\le\xi\circ{\mathbf G}(X,Z)\circ\xi^t\le(1+g_0)|\xi|^2$ for all $\xi\in\mathbb R^3,$ with real $g_0\in[0,+\infty);$
\vspace*{-0.6ex}
\item[(G4)]
$\mbox{det}\,{\mathbf G}(X,Z)=1.$
\end{itemize}
Now, let $X$ be an immersion with Gauss map $N,$ and let ${\mathbf G}(X,Z)$ be a weight matrix as above. We define the components $h_{ij}$ of the weighted first fundamental form of $X$ and the associated line element $ds_g$ as (use the summation convention, and set $u^1\equiv u,$ $u^2\equiv v$)
\begin{equation}\label{2.4}
  h_{ij}:=X_{u^i}\circ{\mathbf G}(X,N)\circ X_{u^j}^t\,,\quad
  ds_g^2=h_{ij}\,du^i du^j\,.
\end{equation}
\subsection{Second step: Definition $\mu$-stability}
We want to prove an area estimate for so-called $\mu$-stable immersions in $\mathbb R^3:$
\begin{definition}
The immersion $X$ is called $\mu$-stable with real $\mu>0$ and a function $q\in C^1(B,\mathbb R)$ iff
\begin{equation}\label{2.5}
  \int\hspace{-0.25cm}\int\limits_{\hspace{-0.3cm}B}\nabla_{ds_g^2}(\varphi,\varphi)W\,dudv
  \ge\mu\int\hspace{-0.25cm}\int\limits_{\hspace{-0.3cm}B}(q-K)W\varphi^2\,dudv
  \quad\mbox{for all}\ \varphi\in C_0^\infty(B,\mathbb R)
\end{equation}
with a weighted metric $ds_g^2$ from (\ref{2.4}), where $q-K\ge 0$ in $B$ with the Gaussian curvature $K$ of the surface, and where
\begin{equation}\label{2.6}
  \nabla_{ds_g^2}(\varphi,\psi):=h^{ij}\,\varphi_{u^i}\psi_{u^j}\,,\quad
  h_{ij}h^{jk}=\delta_i^k\,,
\end{equation}
is the Beltrami operator w.r.t. $ds_g^2,$ $\delta_i^k$ is the Kronecker symbol.
\end{definition}
\begin{examples}\quad\\[1ex]
In the calculus of variations we are faced with ``weighted'' and ``unweighted'' problems:
\vspace*{-0.6ex}
\begin{itemize}
\item[1.]
A conformally parametrized surface of constant mean curvature $h_0\in\mathbb R$ (as a critical point of the area functional with a suitable volume constraint) is stable iff
\begin{equation}\label{2.7}
  \int\hspace{-0.25cm}\int\limits_{\hspace{-0.3cm}B}|\nabla\varphi|^2\,dudv
  \ge 2\int\hspace{-0.25cm}\int\limits_{\hspace{-0.3cm}B}(2h_0^2-K)W\varphi^2\,dudv
  \quad\mbox{for all}\ \varphi\in C_0^\infty(B,\mathbb R),
\end{equation}
that is, it is $\mu$-stable with $\mu=2$ and $q\equiv 2h_0^2.$ We have $q\equiv 0$ for minimal surfaces. Furthermore, in this case ${\mathbf G}(X,Z)\equiv\mathbb E^3$ with the three-dimensional unit matrix $\mathbb E^3\subset\mathbb R^{3\times 3},$ such that $\nabla_{ds^2}(\varphi,\varphi)=\frac{1}{W}\,|\nabla\varphi|^2=\frac{1}{W}\,(\varphi_u^2+\varphi_v^2)$ using conformal parameters, and where $ds^2$ stands for the non-weighted line element.
\vspace*{-0.6ex}
\item[2.]
Critical points $X$ of variational problems
\begin{equation}\label{2.8}
  \int\hspace*{-0.25cm}\int\limits_{\hspace*{-0.3cm}B}F(X,X_u\times X_v)\,dudv
  \longrightarrow\mbox{extr!}
\end{equation}
are immersions of mean-curvature type, that is, they solve
\begin{equation}\label{2.9}
  \nabla_{ds_g^2}(X,N)
  =-2H_G(X,N)
  =-\,\frac{\mbox{trace}\,{\mathbf F}_{XZ}(X,N)}{\sqrt{\mbox{det}\,{\mathbf F}_{ZZ}(X,N)}}\,,
\end{equation}
where ${\mathbf F}_{XZ}=(F_{x^iz^j})_{i,j=1,2,3}\in\mathbb R^{3\times 3}$ etc., with the weighted mean curvature $H_G(X,Z)$ w.r.t.
\begin{equation}\label{2.10}
  {\mathbf G}(X,Z)
  =\left(
     \frac{{\mathbf F}_{ZZ}(X,Z)}{\sqrt{\mbox{det}\,{\mathbf F}_{ZZ}(X,Z)}}
     +(z^iz^j)_{i,j=1,2,3}
   \right)^{-1}\,,
\end{equation}
and $ds_g^2$ chosen as in (\ref{2.4}). This weight matrix was first introduced in \cite{Sauvigny_01}. For example, assume that the integrand in (\ref{2.8}) has the form $F=F(Z).$ If for a critical point (a so-called G-minimal surface, $H_G(X,N)\equiv 0$) the second variation is non-negative, then it can be shown (see \cite{Froehlich_01})
\begin{equation}\label{2.11}
  \int\hspace*{-0.25cm}\int\limits_{\hspace*{-0.3cm}B}\nabla_{ds_g^2}(\varphi,\varphi)W\,dudv
  \ge\mu\int\hspace*{-0.25cm}\int\limits_{\hspace*{-0.3cm}B}(-K)W\varphi^2\,dudv
  \quad\mbox{for all}\ \varphi\in C_0^\infty(B,\mathbb R),
\end{equation}
that is, a critical point is $\mu$-stable with $q\equiv 0$ and a suitable $\mu>0.$
\end{itemize}
Independent of the theory of the second variation, various stability criteria were developed by analysing spherical properties of the immersions:
\begin{itemize}
\item[3.]
For example, stability for minimal surfaces \cite{Barbosa_doCarmo_01}, for surfaces of prescribed constant mean curvature \cite{Ruchert_01}, \cite{Froehlich_02}, for F-minimal surfaces \cite{Clarenz_01}, or for weighted minimal surfaces \cite{Froehlich_01}.
\vspace*{-0.6ex}
\item[4.]
In \cite{Barbosa_doCarmo_02} the reader can find stability criteria for minimal surfaces in the three-sphere $S^3,$ in the hyperbolic space $H^3,$ and in the Euclidean space $\mathbb R^n.$ This last result was improved for minimal graphs with flat normal bundle in \cite{Froehlich_05}.
\end{itemize}
\end{examples}
\subsection{An estimate for the area growth of geodesic discs}
Using methods which go back to \cite{Gulliver_01} and \cite{Sauvigny_01} we prove the following area bound:
\begin{theorem}
Let the immersion $X$ be $\mu$-stable in the sense of (\ref{2.5}), such that
\begin{equation}\label{2.12}
  \mu>\frac{1+g_0}{2}\quad\mbox{and}\quad q\ge 0\ \mbox{in}\ B.
\end{equation}
Let it represent a geodesic disc ${\mathfrak B}_r(X_0)$ of radius $r>0$ and center $X_0\in\mathbb R^3.$ Then
\begin{equation}\label{2.13}
  {\mathcal A}[X]
  \le\frac{2\pi\mu}{2\mu-(1+g_0)}\,r^2
\end{equation}
for the area ${\mathcal A}[X]$ of the immersion.
\end{theorem}
\begin{remarks}
\begin{itemize}
\item[1.]
The assumption $q\ge 0$ in $B$ is needed in the estimate (\ref{2.21}). If it is not fulfilled, we could proceed in (\ref{2.21}) with $q^-(u,v):=\mbox{\rm min}\,\big\{q(u,v),0\big\}.$ It would follow
\begin{equation}\label{2.14}
  {\mathcal A}[X]
  \le\frac{2\pi\mu}{2\mu-(1+g_0)}\,r^2
     -\frac{\mu r^2}{2\mu-(1+g_0)}
      \int\hspace{-0.25cm}\int\limits_{\hspace{-0.3cm}B}q^-(u,v)W(u,v)\,dudv.
\end{equation}
\item[2.]
The proof of the theorem uses intrinsic methods. Therefore, it could be extended to $\mu$-stable immersions in Euclidean spaces $\mathbb R^n$ for $n\ge 3.$ But up to now we are not able to transform critical points of general elliptic variational problems in spaces of higher codimension into a weighted form as given in (\ref{2.9}), (\ref{2.10}) (see also the remarks in section 4.1). 
\item[3.]
The smallest value for $\mu,$ such that a growth estimate of this form is true, is not known; see the discussion in \cite{Colbrie_Schoen_01}.
\end{itemize}
\end{remarks}
\noindent
For the proof we need the following result (see \cite{Froehlich_01}):
\begin{lemma}
Let the immersion $X$ be given. We denote by $ds^2$ its non-weighted line element, and by $ds_g^2$ its weighted element w.r.t. a weight matrix ${\mathbf G}(X,Z).$ Then there hold
\begin{equation}\label{2.15}
\begin{array}{l}
  \displaystyle
  (1+g_0)^{-1}\int\hspace{-0.25cm}\int\limits_{\hspace{-0.3cm}B}\nabla_{ds^2}(\varphi,\varphi)W\,dudv
  \,\le\,\int\hspace{-0.25cm}\int\limits_{\hspace{-0.3cm}B}\nabla_{ds_g^2}(\varphi,\varphi)W\,dudv \\[0.7cm]
  \hspace*{2cm}\displaystyle
  \le\,(1+g_0)\int\hspace{-0.25cm}\int\limits_{\hspace{-0.3cm}B}\nabla_{ds^2}(\varphi,\varphi)W\,dudv
\end{array}
\end{equation}
for all $\varphi\in C_0^1(B,\mathbb R)$ with the Beltrami operators $\nabla_{ds^2}$ and $\nabla_{ds_g^2}$ from (\ref{2.6}).
\end{lemma}
\def\proofname{Proof of the Theorem}
\begin{proof}
\begin{itemize}
\item[1.]
Due to the Lemma, the $\mu$-stability (\ref{2.5}) yields
\begin{equation}\label{2.16}
  \int\hspace{-0.25cm}\int\limits_{\hspace{-0.3cm}B}\nabla_{ds^2}(\varphi,\varphi)\,W\,dudv
  \ge\frac{\mu}{1+g_0}\int\hspace{-0.25cm}\int\limits_{\hspace{-0.3cm}B}(q-K)W\varphi^2\,dudv.
\end{equation}
\item[2.]
Introduce geodesic polar coordinates $(\varrho,\varphi)\in[0,r]\times[0,2\pi].$ For curves $\varrho=\mbox{const}$ on the surface, the integral formula of Bonnet and Gauss reads
\begin{equation}\label{2.17}
  \int\limits_0^{2\pi}\kappa_g(\varrho,\varphi)\sqrt{P(\varrho,\varphi)}\,d\varphi
  +\int\limits_0^\varrho\int\limits_0^{2\pi}K(\tau,\varphi)\sqrt{P(\tau,\varphi)}\,d\tau d\varphi
  =2\pi
\end{equation}
with the geodesic curvature $\kappa_g.$ For the area element $P$ there hold $P(\varrho,\varphi)>0$ for all $(0,r]\times[0,2\pi),$ as well as
\begin{equation}\label{2.18}
  \lim_{\varrho\to 0_+}P(\varrho,\varphi)=0,\quad
  \lim_{\varrho\to 0_+}\frac{\partial}{\partial\varrho}\,\sqrt{P(\varrho,\varphi)}=1
  \quad\mbox{for all}\ \varphi\in[0,2\pi).
\end{equation}
Following \cite{Blaschke_Leichtweiss_01}, \S 81, for such curves it holds $\kappa_g\sqrt{P}=\frac{\partial}{\partial\varrho}\,\sqrt{P}$ for $(\varrho,\varphi)\in(0,r]\times[0,2\pi),$ thus
\begin{equation}\label{2.19}
\begin{array}{lll}
  \displaystyle
  \frac{\partial}{\partial\varrho}\,
  \int\limits_0^{2\pi}\sqrt{P(\varrho,\varphi)}\,d\varphi\negthickspace
  &  =  & \displaystyle\negthickspace
          \int\limits_0^{2\pi}\kappa_g(\varrho,\varphi)\sqrt{P(\varrho,\varphi)}\,d\varphi \\[0.7cm]
  &  =  & \displaystyle\negthickspace
          2\pi-\int\limits_0^\varrho\int\limits_0^{2\pi}K(\tau,\varphi)\sqrt{P(\tau,\varphi)}\,d\tau d\varphi.
\end{array}
\end{equation}
\item[3.]
Define the function $L(\varrho):=\int\limits_0^{2\pi}\sqrt{P(\varrho,\varphi)}\,d\varphi,$ $0<\varrho\le r,$ with the derivatives
\begin{equation}\label{2.20}
  L'(\varrho)
  =2\pi
   -\int\limits_0^\varrho\int\limits_0^{2\pi}
     K(\tau,\varphi)\sqrt{P(\tau,\varphi)}\,d\tau d\varphi,\quad
  L''(\varrho)
  =-\int\limits_0^{2\pi}K(\varrho,\varphi)\sqrt{P(\varrho,\varphi)}\,d\varphi.
\end{equation}
\item[4.]
Consider the test function $\Phi(\varrho):=1-\frac{\varrho}{r},$ $0<\varrho\le r.$ It holds $\nabla_{ds_P^2}(\Phi,\Phi)=\Phi'(\varrho)^2$ with the line element $ds_P^2$ w.r.t. the geodesic polar coordinates. Using $q\ge 0$ we estimate as follows:
\begin{equation}\label{2.21}
\begin{array}{l}
  \displaystyle
  \int\limits_0^r\Phi'(\varrho)^2L(\varrho)\,d\varrho
  \,=\,\int\limits_0^r\int\limits_0^{2\pi}
       \Phi'(\varrho)^2\sqrt{P(\varrho,\varphi)}\,d\varrho d\varphi \\[0.7cm]
  \hspace*{1.6cm}\displaystyle
  \ge\,\frac{\mu}{1+g_0}
       \int\limits_0^r\int\limits_0^{2\pi}
       \Big\{q(\varrho,\varphi)-K(\varrho,\varphi)\Big\}\,\Phi(\varrho)^2\sqrt{P(\varrho,\varphi)}
       \,d\varrho d\varphi\\[0.7cm]
  \hspace*{1.6cm}\displaystyle
  =\,\frac{\mu}{1+g_0}
     \int\limits_0^r\int\limits_0^{2\pi}
     q(\varrho,\varphi)\Phi(\varrho)^2\sqrt{P(\varrho,\varphi)}\,d\varrho d\varphi
     +\frac{\mu}{1+g_0}\int\limits_0^rL''(\varrho)\Phi(\varrho)^2\,d\varrho \\[0.7cm]
  \hspace*{1.6cm}\displaystyle
  \ge\,\frac{\mu}{1+g_0}\int\limits_0^rL''(\varrho)\Phi(\varrho)^2\,d\varrho.
\end{array}
\end{equation}
\item[5.]
Together with (\ref{2.18}), integration by parts yields
\begin{equation}\label{2.22}
\begin{array}{lll}
  \displaystyle
  \int\limits_0^r L''(\varrho)\Phi(\varrho)^2\,d\varrho\negthickspace
  & = & \negthickspace\displaystyle
        L'(\varrho)\Phi(\varrho)\,\Big|_{\varrho=0_+}^{\varrho=r}
         -2\int\limits_0^rL'(\varrho)\Phi(\varrho)\Phi'(\varrho)\,d\varrho \\[0.8cm]
  & = & \negthickspace\displaystyle
        -\,2\pi
        +2\int\limits_0^rL(\varrho)\Phi'(\varrho)^2\,d\varrho.
\end{array}
\end{equation}
Thus, $\int\limits_0^rL(\varrho)\Phi'(\varrho)^2\,d\varrho\le\frac{2\pi\mu}{2\mu-(1+g_0)}\,,$ and the statement follows with $\Phi'^2=\frac{1}{r^2}.$
\end{itemize}
\end{proof}
\begin{example}
Let the immersion $X$ with prescribed constant mean curvature $h_0$ be $\mu$-stable with real $\mu>\frac{1}{2}$ and $q\equiv 2h_0^2$ (compare with (\ref{2.7})). Furthermore, let it represent a geodesic disc ${\mathfrak B}_r(X_0)$ with geodesic radius $r>0$ and center $X_0.$ Then it holds
\begin{equation}\label{2.23}
  {\mathcal A}[X]\le\frac{2\pi\mu}{2\mu-1}\,r^2
\end{equation}
due to ${\mathbf G}(X,Z)\equiv\mathbb E^3,$ that is, $g_0=0.$
\end{example}
\subsection{Remark: Area bounds for minimizers via outer balls}
In \cite{Winklmann_01} we find area bounds in terms of outer balls enclosing embedded minimizers for the general variational problem (\ref{2.8}). Namely, denote by $\nu\colon{\mathcal M}\to S^2$ its unit normal. Intersect the surface with the closed ball $K_\varrho(X_0)$ of radius $\varrho>0$ and center $X_0\in{\mathcal M}.$ Assume that ${\mathcal M}\cap K_\varrho(X_0)$ is simply connected. The greater of the two ``caps'' of the boundary $\partial K_\varrho(X_0),$ which are generated by this intersection, is denoted by ${\mathcal K}.$ Now, assume that
\begin{equation}\label{2.24}
  m_1|Z|\le\widetilde F(X,Z)\le m_2|Z|
  \quad\mbox{for all}\ (X,Z)\in\mathbb R^3\times\mathbb R^3\setminus\{0\}
\end{equation}
for the composition $\widetilde F=F\circ X,$ where $0<m_1\le m_2<+\infty.$ Due to the minimality of $X$ we estimate
\begin{equation}\label{2.25}
\begin{array}{lll}
  \displaystyle
  {\mathcal A}[{\mathcal M}\cap K_\varrho(X_0)]\negthickspace
  &  =  & \negthickspace\displaystyle
          \int\limits_{{\mathcal M}\cap K_\varrho(X_0)}d{\mathcal M}
          \,\le\,\frac{1}{m_1}
                 \int\limits_{{\mathcal M}\cap K_\varrho(X_0)}\widetilde F(X,\nu)\,d{\mathcal M} \\[0.8cm]
  & \le & \negthickspace\displaystyle
          \frac{1}{m_1}\int\limits_{\mathcal K}\widetilde F(X,\nu)\,d{\mathcal K}
          \,\le\,\frac{m_2}{m_1}\int\limits_{\mathcal K}d{\mathcal K}
          \,<\,\frac{4m_2\pi}{m_1}\,\varrho^2\,.
\end{array}
\end{equation}
\subsection{An estimate in terms of the curvatura integra}
The proof of our theorem allows the next result (see also \cite{Sauvigny_01}):
\begin{proposition}
Let the immersion $X$ represent a geodesic disc ${\mathfrak B}_r(X_0)$ of radius $r>0$ and center $X_0.$ Let its Gaussian curvature satisfy
\begin{equation}\label{2.26}
  K(\varrho,\varphi)\le K_0\quad\mbox{for all}\ (\varrho,\varphi)\in[0,r]\times[0,2\pi]
\end{equation}
with a real constant $K_0\in[0,+\infty).$ Then it holds
\begin{equation}\label{2.27}
  {\mathcal A}[X]
  \le r^2
      \left\{
        \pi
        +\frac{1}{2}
         \int\limits_0^r\int\limits_0^{2\pi}
         \Big\{K_0-K(\varrho,\varphi)\Big\}\sqrt{P(\varrho,\varphi)}\,d\varrho d\varphi
      \right\}.
\end{equation}
\end{proposition}
\begin{example}
For minimal surfaces we have $K_0=0.$
\end{example}
\def\proofname{Proof of the Proposition}
\begin{proof}
For curves $\varrho=\mbox{const}$ we conclude from the second line in (\ref{2.19})
\begin{equation}\label{2.28}
  \frac{\partial}{\partial\varrho}
  \int\limits_0^{2\pi}\sqrt{P(\varrho,\varphi)}\,d\varphi
  \le 2\pi
      +\int\limits_0^r\int\limits_0^{2\pi}
       \Big\{K_0-K(\tau,\varphi)\Big\}\sqrt{P(\tau,\varphi)}\,d\tau d\varphi.
\end{equation}
A first integrating w.r.t. the radius coordinate $\varrho,$ and then a further integration w.r.t. to $\varrho=0\ldots r$ proves the statement.
\end{proof}
\subsection{An estimate in terms of the boundary curvature}
Given the immersion $X$ with its $C^2$-regular boundary curve. Denote by $\kappa_g$ and $\kappa_n$ its geodesic curvature and normal curvature, resp. It holds $\kappa=\sqrt{\kappa_g^2+\kappa_n^2}\ge|\kappa_g|$ for the non-negative curvature of the boundary, and due to Bonnet-Gau\ss{} we conclude
\begin{equation}\label{2.29}
  \int\hspace{-0.25cm}\int\limits_{\hspace{-0.3cm}B}(-K)W\,dudv
  =\int\limits_{\partial B}\kappa_g(s)\,ds-2\pi
  \le\int\limits_{\partial B}\kappa(s)\,ds-2\pi.
\end{equation}
Inserting into (\ref{2.27}) proves the
\begin{corollary}
Under the above assumptions it holds
\begin{equation}\label{2.30}
  {\mathcal A}[X]\le\frac{K_0}{2}\,{\mathcal A}[X]r^2+\frac{r^2}{2}\int\limits_{\partial B}\kappa(s)\,ds
\end{equation}
with the constant $K_0\in[0,+\infty)$ from (\ref{2.26}). In particular, if $K_0=0$ then it holds
\begin{equation}\label{2.31}
  {\mathcal A}[X]\le\frac{r^2}{2}\int\limits_{\partial B}\kappa(s)\,ds.
\end{equation}
\end{corollary}
\section{Graphs in $\mathbb R^3$}
\setcounter{equation}{0}
\subsection{Introductory remarks}
In this chapter we want to prove an upper area bound for graphs which solve non-homogeneous quasilinear equations. First, let us give some examples:
\begin{itemize}
\item[1.]
Critical points of variational problems
\begin{equation}\label{3.1}
  \int\hspace*{-0.25cm}\int\limits_{\hspace*{-0.3cm}\Omega}F(x,y,\zeta,\zeta_x,\zeta_y)\,dxdy
  \longrightarrow\mbox{extr!}
\end{equation}
have non-homogeneous divergence form (see section 4.2).
\item[2.]
An equation of the form
\begin{equation}\label{3.2}
   A(\zeta_x,\zeta_y)\zeta_{xx}+2B(\zeta_x,\zeta_y)\zeta_{xy}+C(\zeta_x,\zeta_y)\zeta_{yy}=0
\end{equation}
with smooth coefficients can always be transformed into divergence form (see \cite{Bers_01}).
\item[3.]
By introducing a suitable weight matrix, solutions of
\begin{equation}\label{3.3}
  A(x,y,\zeta,\zeta_x,\zeta_y)\zeta_{xx}
  +2B(x,y,\zeta,\zeta_x,\zeta_y)\zeta_{xy}
  +C(x,y,\zeta,\zeta_x,\zeta_y)\zeta_{yy}
  =R(x,y,\zeta,\zeta_x,\zeta_y)
\end{equation}
can be transformed into the Beltrami form (\ref{2.9}) (see i.e. \cite{Sauvigny_01}, \cite{Sauvigny_02} for $R\equiv 0$). This would make the results of chapter 2 applicable to the objects of study in this part (see also the examples discussed in section 2.2).
\end{itemize}
\subsection{An estimate in the general case}
The next result follows ideas from \cite{Finn_01}, where area estimates for homogeneous divergence form equations are established, but where an explicit form as below is not needed. Therefore, we want to demonstrate all the essential steps.
\begin{theorem}
Let $\zeta\in C^2(\Omega,\mathbb R)\cap C^1(\overline\Omega,\mathbb R),$ $\Omega\subset\mathbb R^2$ bounded and simply connected and with $C^1$-regular boundary, solve the elliptic Dirichlet boundary value problem
\begin{equation}\label{3.4}
\begin{array}{rcl}
  \displaystyle
  \frac{d}{d x}\,F_p(x,y,\zeta,\zeta_x,\zeta_y)
  +\frac{d}{dy}\,F_q(x,y,\zeta,\zeta_x,\zeta_y)\negthickspace
  & = & \negthickspace
        R(x,y,\zeta,\zeta_x,\zeta_y)\quad\mbox{in}\ \Omega, \\[0.3cm]
  \zeta(x,y)\negthickspace
  & = & \negthickspace
        \varphi(x,y)\quad\mbox{on}\ \partial\Omega,
\end{array}
\end{equation}
where $\varphi\in C^1(\mathbb R^2,\mathbb R).$ Assume that
\begin{itemize}
\item[(A1)]
for all $(x,y,z,p,q)\in\mathbb R^5$
\begin{equation}\label{3.5}
  F_p(x,y,z,p,q)^2+F_q(x,y,z,p,q)^2\le k_0^2
\end{equation}
with a real constant $k_0\in[0,+\infty);$
\item[(A2)]
with a further real constant $m_1\in(0,+\infty)$
\begin{equation}\label{3.6}
  m_1|\xi|^2
  \le(\xi_1,\xi_2)\circ
     \begin{pmatrix}
       F_{pp}(x,y,z,\tilde p,\tilde q) & F_{pq}(x,y,z,\tilde p,\tilde q) \\[0.1cm]
       F_{qp}(x,y,z,\tilde p,\tilde q) & F_{qq}(x,y,z,\tilde p,\tilde q)
     \end{pmatrix}\circ
     \begin{pmatrix} \xi_1 \\[0.1cm] \xi_2 \end{pmatrix}
\end{equation}
for all $\xi=(\xi_1,\xi_2)\in\mathbb R^2$ and all $\tilde p,\tilde q\in\mathbb R$ such that $\tilde p^2+\tilde q^2\le 1;$
\item[(A3)]
finally
\begin{equation}\label{3.7}
  F_p(x,y,z,0,0)=0,\quad F_q(x,y,z,0,0)=0.
\end{equation}
\end{itemize}
Then it holds
\begin{equation}\label{3.8}
  {\mathcal A}[\zeta]
  \le\left(1+\frac{\|\zeta\|_{0,\Omega}\|R\|_{0,\Omega}}{m_1}\right){\mathcal A}\,[\Omega]
    +\frac{\|\zeta\|_{0,\partial\Omega}k_0}{m_1}\,{\mathcal L}[\partial\Omega]
\end{equation}
with the area ${\mathcal A}\,[\Omega]$ of $\Omega\subset\mathbb R^2$ and the length ${\mathcal L}[\partial\Omega]$ of its boundary curve $\partial\Omega,$ and the usual Schauder norms $\|\cdot\|_{0,\Omega}$ etc.
\end{theorem}
\def\proofname{Proof of the Theorem}
\begin{proof}
\begin{itemize}
\item[1.]
Consider the function
\begin{equation}\label{3.9}
  \mu(t):=pF_p(x,y,z,tp,tq)+qF_q(x,y,z,tp,tq),\quad t\in[0,1].
\end{equation}
Assumption (A3) implies $\mu(0)=0$ and $\mu(1)=pF_p(x,y,z,p,q)+qF_q(x,y,z,p,q).$
\item[2.]
For real $t\in[0,1]$ we introduce a real number $m_1^*=m_1^*(t)\in(0,+\infty)$ such that
\begin{equation}\label{3.10}
  m_1^*(t)|\xi|^2
  \le(\xi_1,\xi_2)\circ
     \left(\begin{array}{cc}
       F_{pp}(x,y,z,tp,tq) & F_{pq}(x,y,z,tp,tq) \\[0.1cm]
       F_{pq}(x,y,z,tp,tq) & F_{qq}(x,y,z,tp,tq)
     \end{array}\right)\circ
     \left(\begin{array}{c}
       \xi_1 \\[0.1cm] \xi_2
     \end{array}\right)
\end{equation}
for all $\xi=(\xi_1,\xi_2)\in\mathbb R^2.$ Namely, due to (A2) we demand
\begin{itemize}
\item[$(\alpha)$]
if $p^2+q^2\ge 1,$ then $m_1^*(t)\ge m_1$ for $\displaystyle t\le\frac{1}{\sqrt{p^2+q^2}}\,;$
\item[$(\beta)$]
if $p^2+q^2\le 1,$ then $m_1^*(t)\ge m_1$ for $t\le 1.$
\end{itemize}
Note that in both cases $t^2p^2+t^2q^2\le 1.$
\item[3.]
Differentiating $\mu=\mu(t)$ yields
\begin{equation}\label{3.11}
  \mu'(t)=F_{pp}(x,y,z,tp,tq)p^2+2F_{pq}(x,y,z,tp,tq)pq+F_{qq}(x,y,z,tp,tq)q^2\,,
\end{equation}
and by definition of $m_1^*(t)$ we have $m_1^*(t)(p^2+q^2)\le\mu'(t).$ It follows that
\begin{equation}\label{3.12}
  \mu(1)=\int\limits_0^1\mu'(t)\,dt\ge(p^2+q^2)\int\limits_0^1m_1^*(t)\,dt.
\end{equation}
\item[4.]
Now, note that
\begin{itemize}
\item[$(\gamma)$]
if $p^2+q^2\ge 1,$ then due to ($\alpha$)
\begin{equation}\label{3.13}
  \mu(1)
  \ge(p^2+q^2)\!\!\int\limits_0^{(p^2+q^2)^{-\frac{1}{2}}}m_1^*(t)\,dt
     \ge(p^2+q^2)\!\!\int\limits_0^{(p^2+q^2)^{-\frac{1}{2}}}m_1\,dt
     =m_1\sqrt{p^2+q^2}\,;
\end{equation}
\item[$(\delta)$]
if $p^2+q^2\le 1,$ then due to ($\beta$)
\begin{equation}\label{3.14}
  \mu(1)
  \ge(p^2+q^2)\int\limits_0^1m_1\,dt
  =m_1(p^2+q^2).
\end{equation}
\end{itemize}
Summarising we arrive at (cp. \cite{Finn_01}, Lemma 4)
\begin{equation}\label{3.15}
  pF_p(x,y,z,p,q)+qF_q(x,y,z,p,q)
  \ge\left\{
       \begin{array}{ll}
         m_1(p^2+q^2),            & \mbox{if}\ p^2+q^2\le 1 \\[0.2cm]
         m_1\,\sqrt{p^2+q^2}\,,   & \mbox{if}\ p^2+q^2\ge 1
       \end{array}
     \right.\hspace*{-0.1cm}.
\end{equation}
\item[5.]
Making use of the divergence structur of our Dirichlet problem we infer
\begin{equation}\label{3.16}
  \mbox{div}\,(\zeta F_p,\zeta F_q)
  =pF_p+qF_q+\zeta\left(\frac{d}{dx}\,F_p+\frac{d}{dy}\,F_q\right)
  =pF_p+qF_q+\zeta R\,,
\end{equation}
and integration by parts yields (cp. \cite{Finn_01}, Lemma 5)
\begin{equation}\label{3.17}
\begin{array}{l}
  \displaystyle
  \int\hspace{-0.25cm}\int\limits_{\hspace{-0.3cm}\Omega}(pF_p+qF_q)\,dxdy
  \,=\,\int\hspace{-0.25cm}\int\limits_{\hspace{-0.3cm}\Omega}
       \mbox{div}\,(\zeta F_p,\zeta F_q)\,dxdy
       -\int\hspace{-0.25cm}\int\limits_{\hspace{-0.3cm}\Omega}\zeta R\,dxdy \\[0.8cm]
  \hspace*{1.2cm}\displaystyle
  =\,\int\limits_{\partial\Omega}\zeta(F_p,F_q)\cdot\nu^t\,ds
     -\int\hspace{-0.25cm}\int\limits_{\hspace{-0.3cm}\Omega}\zeta R\,dxdy \\[0.8cm]
  \hspace*{1.2cm}\displaystyle
  \le\,\|\zeta\|_{0,\partial\Omega}
       \int\limits_{\partial\Omega}\sqrt{F_p^2+F_q^2}\,ds
       +\|\zeta\|_{0,\Omega}\|R\|_{0,\Omega}\,{\mathcal A}\,[\Omega] \\[0.8cm]
  \hspace*{1.2cm}\displaystyle
  \le\,\|\zeta\|_{0,\partial\Omega}k_0\,{\mathcal L}[\partial\Omega]
       +\|\zeta\|_{0,\Omega}\|R\|_{0,\Omega}\,{\mathcal A}\,[\Omega]
\end{array}
\end{equation}
with $\nu=\nu(s)$ normal to the boundary $\partial\Omega\subset\mathbb R^2.$
\item[6.]
Taking $\sqrt{1+p^2+q^2}\le 1+\sqrt{p^2+q^2}$ for $p^2+q^2\ge 1,$ and $\sqrt{1+p^2+q^2}\le 1+p^2+q^2$ for $p^2+q^2\le 1$ into account, we calculate for $p^2+q^2\le 1$ and $p^2+q^2\ge 1$ (cp. \cite{Finn_01}, Proof of Theorem III)
\begin{equation}\label{3.18}
  \int\hspace{-0.25cm}\int\limits_{\hspace{-0.3cm}\Omega}\sqrt{1+p^2+q^2}\,dxdy
  \le{\mathcal A}\,[\Omega]
     +\frac{1}{m_1}\,
      \int\hspace{-0.25cm}\int\limits_{\hspace{-0.3cm}\Omega}(pF_p+qF_q)\,dxdy.
\end{equation}
The statement follows.
\end{itemize}
\end{proof}
\begin{remark}
Various variations of the proof are possible: For example, we could alter (\ref{3.17}) to obtain
  $$\int\hspace{-0.25cm}\int\limits_{\hspace{-0.3cm}\Omega}(pF_p+qF_q)\,dxdy
    \le\|\zeta\|_{0,\partial\Omega}k_0\,{\mathcal L}[\partial\Omega]
       +\|\zeta\|_{0,\Omega}\|R\|_{L^1(\Omega)}$$
with the $L^1$-norm on $\Omega.$ Then (\ref{3.8}) would change according to this new estimate.
\end{remark}
\subsection{Homogeneous divergence equations}
\begin{corollary}
In the homogeneous case $R\equiv 0$ we conclude from the Theorem
\begin{equation}\label{3.19}
  {\mathcal A}[\zeta]
  \le{\mathcal A}[\Omega]
     +\frac{\|\zeta\|_{0,\partial\Omega}k_0}{m_1}\,{\mathcal L}[\partial\Omega].
\end{equation}
\end{corollary}
\noindent
For example, let us consider minimal graphs with $F(p,q)=\sqrt{1+p^2+q^2}$ such that
\begin{equation}\label{3.20}
  F_p(p,q)=\frac{p}{\sqrt{1+p^2+q^2}}\,,\quad
  F_q(p,q)=\frac{q}{\sqrt{1+p^2+q^2}}\,.
\end{equation}
For (A1), we calculate $F_p^2+F_q^2=\frac{p^2+q^2}{1+p^2+q^2}\le 1=:k_0.$ Furthermore, we set $m_1:=\frac{1}{\sqrt{8}}$ due to
\begin{equation}\label{3.21}
  \frac{1}{\sqrt{8}}\cdot|\xi|^2
  \le(\xi_1,\xi_2)\circ
  \begin{pmatrix}
    \displaystyle
    \frac{1+q^2}{(1+p^2+q^2)^\frac{3}{2}}
      & \displaystyle
        -\,\frac{pq}{(1+p^2+q^2)^\frac{3}{2}}   \\[0.6cm]
    \displaystyle
    -\,\frac{pq}{(1+p^2+q^2)^\frac{3}{2}}
      & \displaystyle
        \frac{1+p^2}{(1+p^2+q^2)^\frac{3}{2}}
  \end{pmatrix}\circ
  \begin{pmatrix}
    \xi_1 \\[0.3cm] \xi_2
  \end{pmatrix}
\end{equation}
as well as $\lambda_1=\frac{1}{(1+p^2+q^2)^\frac{3}{2}}\ge\frac{1}{\sqrt{8}},$ $p^2+q^2\le 1,$ for the ``restricted'' smallest eigenvalue from (A2).
\begin{corollary}
For minimal graphs it holds
\begin{equation}\label{3.22}
  {\mathcal A}[\zeta]
  \le{\mathcal A}[\Omega]
     +\sqrt{8}\,\|\zeta\|_{0,\partial\Omega}{\mathcal L}[\partial\Omega].
\end{equation}
\end{corollary}
\noindent
For the inhomogeneous divergence equation
\begin{equation}\label{3.23}
  \mbox{div}\,\frac{(p,q)}{\sqrt{1+p^2+q^2}}=2H(x,y,z)
\end{equation}
with prescribed mean curvature $H$ such that $h_0=\|H\|_{0,\Omega},$ we conclude
\begin{corollary}
In this case of prescribed mean curvature it holds
\begin{equation}\label{3.24}
  {\mathcal A}[\zeta]
  \le\Big\{1+2\sqrt{8}\,h_0\|\zeta\|_{0,\Omega}\Big\}\,{\mathcal A}[\Omega]
     +\sqrt{8}\,\|\zeta\|_{0,\partial\Omega}{\mathcal L}[\partial\Omega].
\end{equation}
\end{corollary}
\begin{remark}
This result is not sharp. From the estimate of the next chapter we will conclude
\begin{equation}\label{3.25}
  {\mathcal A}[\zeta]
  \le\Big\{1+2h_0\|\zeta\|_{0,\Omega}\Big\}\,\mbox{\rm Area}\,[\Omega]
     +\|\zeta\|_{0,\partial\Omega}{\mathcal L}[\partial\Omega].
\end{equation}
\end{remark}
\subsection{An interior estimate}
The next result is motivated from \cite{Colding_01} where sharp bounds for mean-curvature-graphs are proved.
\begin{proposition}
For real $\nu>0$ we define the interior set
\begin{equation}\label{3.26}
  \Omega_\nu:=\Big\{(x,y)\in\Omega\,:\,\mbox{\rm dist}\,\big((x,y),\partial\Omega\big)>\nu\Big\}.
\end{equation}
Then, under the conditions of the above Theorem and the additional assumption
\begin{equation}\label{3.27}
  pF_p(x,y,z,p,q)+qF_q(x,y,z,p,q)\ge 0\quad\mbox{for all}\ (x,y,z,p,q)\in\mathbb R^5
\end{equation}
(compare with (\ref{3.20})) it holds
\begin{equation}\label{3.28}
  \int\hspace{-0.25cm}\int\limits_{\hspace{-0.32cm}\Omega_\nu}\sqrt{1+\zeta_x^2+\zeta_y^2}\,dxdy
  \le{\mathcal A}[\Omega]
     +\frac{1}{m_1}
      \left(
        \frac{2k_0}{\nu}+\|R\|_{0,\Omega}
     \right)\|\zeta\|_{0,\Omega}\,{\mathcal A}[\Omega].
\end{equation}
\end{proposition}
\def\proofname{Proof}
\begin{proof}
Choose a test function $\varphi\in C_0^\infty(\Omega,\mathbb R)$ such that
\begin{equation}\label{3.29}
  \varphi(u,v)=1\quad\mbox{in}\ \Omega_\nu\,,\quad
  |\nabla\varphi(u,v)|\le\frac{2}{\nu}\quad\mbox{in}\ \Omega.
\end{equation}
We compute $\mbox{div}(\varphi\zeta F_p,\varphi\zeta F_q)=\zeta\nabla\varphi\cdot(F_p,F_q)^t+\varphi(pF_p+qF_q)+\varphi\zeta R.$ Integrating the divergence term would gives no contribution due to $\varphi=0$ on $\partial\Omega.$ Therefore,
\begin{equation}\label{3.30}
\begin{array}{lll}
  \displaystyle
  \int\hspace{-0.25cm}\int\limits_{\hspace{-0.32cm}\Omega_\nu}(pF_p+qF_q)\,dxdy\negthickspace
  & \le  & \negthickspace\displaystyle
          \int\hspace{-0.25cm}\int\limits_{\hspace{-0.3cm}\Omega}\varphi(pF_p+qF_q)\,dxdy \\[0.8cm]
  &  =  & \negthickspace\displaystyle
          -\,\int\hspace{-0.25cm}\int\limits_{\hspace{-0.3cm}\Omega}
             \zeta\nabla\varphi\cdot(F_p,F_q)^t\,dxdy
          -\int\hspace{-0.25cm}\int\limits_{\hspace{-0.3cm}\Omega}
           \varphi\zeta R\,dxdy \\[0.8cm]
  & \le & \negthickspace\displaystyle
          \left(
            \frac{2k_0}{\nu}+\|R\|_{0,\Omega}
          \right)\|\zeta\|_{0,\Omega}\,{\mathcal A}[\Omega].
\end{array}
\end{equation}
We proceed as in point 6, i.e. (\ref{3.15}) and (\ref{3.18}), of the proof of our theorem.
\end{proof}
\section{Fermat-type graphs in $\mathbb R^n$}
\setcounter{equation}{0}
\subsection{Introductory remarks}
In this final chapter we establish an area bound for graphs of Fermat-type in divergence form which are critical for the variational problem ($X=(x,y,\zeta_1,\ldots,\zeta_{n-2})$)
\begin{equation}\label{4.1}
  \int\hspace*{-0.25cm}\int\limits_{\hspace*{-0.3cm}\Omega}\Gamma(X)W\,dxdy
  \longrightarrow\mbox{extr!}
\end{equation}
\begin{itemize}
\item[1.]
For $\Gamma(X)\equiv 1$ we have the usual area functional.
\vspace*{-0.6ex}
\item[2.]
In contrast to the case of $n=3,$ the following area bounds depend additionally on the derivatives of the graphs on the boundary. For example, the area of the conformally parametrized minimal graph $(z,z^n),$ $z=x+iy\in B$ and $n\in\mathbb N,$ depends on the maximum norm of the mapping (which does not depend on $n$) and the exponent $n.$
\vspace*{-0.6ex}
\item[3.]
It remains open how to transform the Euler-Lagrange system of (\ref{4.1}) into a Beltrami form by means of a suitable weight matrix (see the remarks in section 3.1).
\end{itemize}
\subsection{The Euler-Lagrange equations}
Let us start with the general functional
\begin{equation}\label{4.2}
  {\mathcal F}[\zeta_1,\ldots,\zeta_{n-2}]
  =\int\hspace{-0.25cm}\int\limits_{\hspace{-0.3cm}\Omega}
   F(x,y,\zeta_1,\ldots,\zeta_{n-2},\nabla\zeta_1,\ldots,\nabla\zeta_{n-2})\,dxdy.
\end{equation}
We set $\zeta=(\zeta_1,\ldots,\zeta_{n-2}),$ $p_\sigma=\zeta_{\sigma,x},$ $q_\sigma=\zeta_{\sigma,y}$ etc.
\begin{proposition}
The $n-2$ Euler-Lagrange equations of ${\mathcal F}[\zeta_1,\ldots,\zeta_{n-2}]$ are
\begin{equation}\label{4.3}
  \frac{dF_{p_\sigma}(x,y,\zeta,\nabla\zeta)}{dx}
  +\frac{dF_{q_\sigma}(x,y,\zeta,\nabla\zeta)}{dy}
  =F_{z_\sigma}(x,y,\zeta,\nabla\zeta)
  \quad\mbox{for}\ \sigma=1,\ldots,n-2.
\end{equation}
\end{proposition}
\begin{corollary}
The non-parametric minimal surface system is
\begin{equation}\label{4.5}
  \mbox{\rm div}\,\frac{(p_\sigma,q_\sigma)}{W}
  =-\,\mbox{\rm div}\,
   \frac{\displaystyle
         \left(
           p_\sigma\sum_{\theta=1}^{n-2}q_\theta^2
           -q_\sigma\sum_{\theta=1}^{n-2}p_\theta q_\theta\,,\,
           q_\sigma\sum_{\theta=1}^{n-2}p_\theta^2
           -p_\sigma\sum_{\theta=1}^{n-2}p_\theta q_\theta
         \right)
        }{W}
\end{equation}
for $\sigma=1,\ldots,n-2.$
\end{corollary}
\def\proofname{Proof of the Corollary}
\begin{proof}
From $X_x=(1,0,\zeta_{1,x},\ldots,\zeta_{n-2,x}),$ $X_y=(0,1,\zeta_{1,y},\ldots,\zeta_{n-2,y})$ it follows
\begin{equation}\label{4.6}
  h_{11}=1+\sum_{\sigma=1}^{n-2}\zeta_{\sigma,x}^2=1+p^2\,,\quad
  h_{12}=\sum_{\sigma=1}^{n-2}\zeta_{\sigma,x}\zeta_{\sigma,y}=p\cdot q^t\,,\quad
  h_{22}=1+\sum_{\sigma=1}^{n-2}\zeta_{\sigma,y}^2=1+q^2
\end{equation}
setting $p=(p_1,\ldots,p_{n-2}),$ $q=(q_1,\ldots,q_{n-2}).$ Therefore, we have (let $W=F(p,q)$)
\begin{equation}\label{4.7}
  {\mathcal A}[\zeta]
  =\int\hspace*{-0.25cm}\int\limits_{\hspace*{-0.3cm}\Omega}
   F(p,q)\,dxdy
  \equiv\int\hspace*{-0.25cm}\int\limits_{\hspace*{-0.3cm}\Omega}
        \sqrt{1+p^2+q^2+p^2q^2-(p\cdot q^t)^2}\,dxdy.
\end{equation}
Differentiation shows
\begin{equation}\label{4.8}
  F_{p_\sigma}
  =\frac{p_\sigma+p_\sigma q^2-q_\sigma(p\cdot q^t)}{W}\,,\quad
  F_{q_\sigma}
  =\frac{q_\sigma+q_\sigma p^2-p_\sigma(p\cdot q^t)}{W}\,,
\end{equation}
together with $F_{z_\sigma}(p,q)\equiv 0$ for $\sigma=1,\ldots,n-2.$ The statement follows.
\end{proof}
\begin{remark}
For $n=3,$ the minimal surface system reduces to $\mbox{div}\,\frac{\nabla\zeta}{W}=0$ in $\Omega.$
\end{remark}
\begin{corollary}
Critical points $(x,y,\zeta)$ of Fermat's functional with the integrand
\begin{equation}\label{4.9}
  F(x,y,z,p,q)=\Gamma(x,y,z)\sqrt{1+p^2+q^2+p^2q^2-(p\cdot q^t)^2}
\end{equation}
solve the Euler-Lagrange system
\begin{equation}\label{4.10}
\begin{array}{lll}
  \displaystyle
  \mbox{\rm div}\,\frac{(p_\sigma,q_\sigma)}{W}\negthickspace
  & = & \negthickspace\displaystyle
        2H(X,\widetilde N_\sigma)\,\frac{\sqrt{1+p_\sigma^2+q_\sigma^2}}{W} \\[0.4cm]
  &   & \negthickspace\displaystyle
        +\,\frac{1}{\Gamma W}
           \left\{
             \Big[\,p^2+q^2+p^2q^2-(p\cdot q^t)^2\Big]\,\Gamma_{z_\sigma}
             -\sum_{\omega=1}^{n-2}(p_\sigma p_\omega+q_\sigma q_\omega)\Gamma_{z_\omega}
           \right\} \\[0.8cm]
  &   & \negthickspace\displaystyle
        -\,\frac{1}{\Gamma}\,
           \mbox{\rm div}
           \left(
             \frac{p_\sigma q^2-q_\sigma(p\cdot q^t)}{W}\,\Gamma,
             \frac{q_\sigma p^2-p_\sigma(p\cdot q^t)}{W}\,\Gamma
           \right)
\end{array}
\end{equation}
for $\sigma=1,\ldots,n-2$ with the mean curvature field
\begin{equation}\label{4.11}
  H(X,\widetilde N_\sigma)
  =\frac{\Gamma_X(X)\cdot\widetilde N_\sigma^t}{2\Gamma(X)W}\,,\quad X=(x,y,\zeta),
\end{equation}
w.r.t. to the non-othogonally unit normal field
\begin{equation}\label{4.12}
  \widetilde N_\sigma
  =\frac{1}{\sqrt{1+|\nabla\zeta_\sigma|^2}}\,
   (-\zeta_{\sigma,x},-\zeta_{\sigma,y},0,\ldots,0,1,0,\ldots,0),\quad
  \sigma=1,\ldots,n-2.
\end{equation}
\end{corollary}
\begin{proof}
We compute
\begin{equation}\label{4.13}
  F_{p_\sigma}
  =\frac{p_\sigma+p_\sigma q^2-q_\sigma(p\cdot q^t)}{W}\,\Gamma,\quad
  F_{q_\sigma}
  =\frac{q_\sigma+q_\sigma p^2-p_\sigma(p\cdot q^t)}{W}\,\Gamma,\quad
  F_{z_\sigma}=\Gamma_{z_\sigma}W
\end{equation}
as well as
\begin{equation}\label{4.14}
\begin{array}{lll}
  \displaystyle
  \frac{dF_{p_\sigma}}{dx}\negthickspace
  & = & \negthickspace\displaystyle
        \Gamma\,\frac{d}{dx}\,\frac{p_\sigma}{W}
        +\Gamma\,\frac{d}{dx}\,\frac{p_\sigma q^2-q_\sigma(p\cdot q^t)}{W}
        +\frac{p_\sigma+p_\sigma q^2-q_\sigma(p\cdot q^t)}{W}\,\frac{d\Gamma}{dx}\,, \\[0.6cm]
  \displaystyle
  \frac{dF_{q_\sigma}}{dy}\negthickspace
  & = & \negthickspace\displaystyle
        \Gamma\,\frac{d}{dy}\,\frac{q_\sigma}{W}
        +\Gamma\,\frac{d}{dy}\,\frac{q_\sigma p^2-p_\sigma(p\cdot q^t)}{W}
        +\frac{q_\sigma+q_\sigma p^2-p_\sigma(p\cdot q^t)}{W}\,\frac{d\Gamma}{dy}\,.
\end{array}
\end{equation}
Thus, (\ref{4.3}) takes the form
\begin{equation}\label{4.15}
\begin{array}{lll}
  \displaystyle
  \mbox{div}\,\frac{(p_\sigma,q_\sigma)}{W}\negthickspace
  & = & \negthickspace\displaystyle
        \frac{\Gamma_{z_\sigma}W}{\Gamma}
        -\mbox{div}
         \left(
           \frac{p_\sigma q^2-q_\sigma(p\cdot q^t)}{W}\,,
           \frac{q_\sigma p^2-p_\sigma(p\cdot q^t)}{W}
         \right) \\[0.6cm]
  &   & \negthickspace\displaystyle
        -\,\frac{p_\sigma+p_\sigma q^2-q_\sigma(p\cdot q^t)}{\Gamma W}\,\frac{d\Gamma}{dx}
        -\frac{q_\sigma+q_\sigma p^2-p_\sigma(p\cdot q^t)}{\Gamma W}\,\frac{d\Gamma}{dy}\,.
\end{array}
\end{equation}
Performing the differentiation gives ($\frac{\partial\Gamma}{\partial x}=\Gamma_x+\Gamma_{z_\omega}p_\omega$ etc.)
\begin{equation}\label{4.16}
\begin{array}{lll}
  \displaystyle
  \mbox{div}\,\frac{(p_\sigma,q_\sigma)}{W}\negthickspace
  & = & \negthickspace\displaystyle
        \frac{1}{\Gamma W}\,\Big\{\Gamma_{z_\sigma}-p_\sigma\Gamma_x-q_\sigma\Gamma_y\Big\} \\[0.4cm]
  &   & \negthickspace\displaystyle
        +\,\frac{1}{\Gamma W}
           \left\{
             \Big[\,p^2+q^2+p^2q^2-(p\cdot q^t)^2\Big]\,\Gamma_{z_\sigma}
             -\sum_{\omega=1}^{n-2}(p_\sigma p_\omega+q_\sigma q_\omega)\Gamma_{z_\omega}
           \right\} \\[0.8cm]
  &   & \negthickspace\displaystyle
        -\,\frac{p_\sigma q^2-q_\sigma(p\cdot q^t)}{\Gamma W}\,\frac{d\Gamma}{dx}
        -\frac{q_\sigma p^2-p_\sigma(p\cdot q^t)}{\Gamma W}\,\frac{d\Gamma}{dy} \\[0.6cm]
  &   & \negthickspace\displaystyle
        -\,\mbox{div}
           \left(
             \frac{p_\sigma q^2-q_\sigma(p\cdot q^t)}{W}\,,
             \frac{q_\sigma p^2-p_\sigma(p\cdot q^t)}{W}
           \right).
\end{array}
\end{equation}
With $\Gamma_X=(\Gamma_x,\Gamma_y,\Gamma_{z_1},\ldots,\Gamma_{z_{n-2}}),$ the first row can be transformed into
\begin{equation}\label{4.17}
  \frac{1}{\Gamma W}\,\Big\{\Gamma_{z_\sigma}-p_\sigma\Gamma_x-q_\sigma\Gamma_y\Big\}
  =\frac{1}{\Gamma W}\,\widetilde N_\sigma\cdot\Gamma_X^t\sqrt{1+|\nabla\zeta_\sigma|^2}
  =2H(X,\widetilde N_\sigma)\,\frac{\sqrt{1+|\nabla\zeta_\sigma|^2}}{W}\,,
\end{equation}
and the statement follows.
\end{proof}
\begin{remark}
For $n=3,$ the Euler-Lagrange system reduces to
\begin{equation}\label{4.18}
  \mbox{div}\,\frac{\nabla\zeta}{W}=2H(X,\widetilde N)\quad\mbox{in}\ \Omega.
\end{equation}
\end{remark}
\subsection{An area estimate}
The main result of this chapter is the following
\begin{theorem}
Let $\zeta\in C^1(\overline\Omega,\mathbb R^{n-2})\cap C^2(\Omega,\mathbb R^{n-2})$ solve (\ref{4.1}) where $\Gamma=\Gamma(x,y)\in C^1(\overline\Omega,\mathbb R)$ such that with real constants $\Gamma_0,$ $\Gamma_1,$ and $\Gamma_2$ it holds
\begin{equation}\label{4.19}
  0<\Gamma_0\le\Gamma(X)\le\Gamma_1<+\infty,\quad
  \Gamma_2:=\|\Gamma\|_{1,\Omega}\,,\quad
  h_0:=\sup_{(X,Z)\in\mathbb R^3\times S^1}|H(X,Z)|.
\end{equation}
We require the smallness condition
\begin{equation}\label{4.20}
  \Lambda:=1-\frac{\sqrt{2}\,(n-2)^2\Gamma_2}{\Gamma_0}\,\max_{\sigma=1,\ldots,n-2}\|\zeta_\sigma\|_{0,\Omega}>0.
\end{equation}
Then it holds
\begin{equation}\label{4.21}
\hspace*{-0.6cm}
\begin{array}{lll}
  \Lambda\cdot{\mathcal A}[\zeta]\negthickspace
  & \le & \negthickspace\displaystyle
          {\mathcal A}[\Omega]
          +(n-2){\mathcal L}[\partial\Omega]\max_{\sigma=1,\ldots,n-2}\|\zeta_\sigma\|_{0,\partial\Omega} \\[0.4cm]
  &     & \negthickspace\displaystyle
          +\,2(n-2)h_0\,{\mathcal A}[\Omega]\max_{\sigma=1,\ldots,n-2}\|\zeta_\sigma\|_{0,\Omega} \\[0.4cm]
  &     & \negthickspace\displaystyle
          +\,(n-2)^2{\mathcal L}[\partial\Omega]
             \max_{\sigma=1,\ldots,n-2}
             \|\zeta_\sigma\|_{0,\partial\Omega}\|D^\top\zeta_\sigma\|_{0,\partial\Omega}
\end{array}
\end{equation}
with the tangential derivative $D^\top\zeta_\sigma=(q_\sigma,-p_\sigma)\cdot v^t,$ $\nu$ unit normal along $\partial\Omega,$ $\sigma=1,\ldots,n-2.$
\end{theorem}
\begin{remarks}\quad
\begin{itemize}
\item[1.]
The third line in  (\ref{4.21}) does not appear if $n=3.$ Furthermore, in this case we set $\Lambda:=1.$ Furthermore, $\Gamma_2=0$ implies $\Lambda=1.$
\vspace*{-0.6ex}
\item[2.]
If we prescribe boundary values $\zeta_{\sigma,R},$ we can replace the tangential derivatives $D^\top\zeta_\sigma$ by the derivatives of $\zeta_{\sigma,R}.$
\end{itemize}
\end{remarks}
\def\proofname{Proof of the Theorem}
\begin{proof}
\begin{itemize}
\item[1.]
We add the $n-2$ identities $\zeta_\sigma\,\mbox{div}\,\frac{\nabla\zeta_\sigma}{W}=\mbox{div}\,\frac{\zeta_\sigma\nabla\zeta_\sigma}{W}-\frac{|\nabla\zeta_\sigma|^2}{W}:$
\begin{equation}\label{4.22}
  \sum_{\sigma=1}^{n-2}\zeta_\sigma\,\mbox{div}\,\frac{\nabla\zeta_\sigma}{W}
  =\sum_{\sigma=1}^{n-2}\mbox{div}\,\frac{\zeta_\sigma\nabla\zeta_\sigma}{W}
   -\sum_{\sigma=1}^{n-2}\frac{p_\sigma^2+q_\sigma^2}{W}
  =\sum_{\sigma=1}^{n-2}\mbox{div}\,\frac{\zeta_\sigma\nabla\zeta_\sigma}{W}
   -\frac{p^2+q^2}{W}\,.
\end{equation}
For the area element we have
\begin{equation}\label{4.23}
  \frac{1}{W}-W
  =\frac{1-\big[1+p^2+q^2+p^2q^2-(p\cdot q^t)^2\big]}{W}
  =-\,\frac{p^2+q^2}{W}-\frac{p^2q^2-(p\cdot q^t)^2}{W}\,,
\end{equation}
therefore,
\begin{equation}\label{4.24}
  W=\frac{1}{W}
    +\sum_{\sigma=1}^{n-2}\mbox{div}\,\frac{\zeta_\sigma\nabla\zeta_\sigma}{W}
    -\sum_{\sigma=1}^{n-2}\zeta_\sigma\,\mbox{div}\,\frac{\nabla\zeta_\sigma}{W}
    +\frac{p^2q^2-(p\cdot q^t)^2}{W}\,.
\end{equation}
\item[2.]
Multiply the Euler-Lagrange equations (\ref{4.10}) by $\zeta_\sigma.$ Summation gives ($\Gamma_{z_\sigma}\equiv 0$!)
\begin{equation}\label{4.25}
\begin{array}{lll}
  \displaystyle
  \sum_{\sigma=1}^{n-2}\zeta_\sigma\,\mbox{div}\,\frac{\nabla\zeta_\sigma}{W}\negthickspace
  & = & \negthickspace\displaystyle
        2\sum_{\sigma=1}^{n-2}
        H(X,\widetilde N_\sigma)\zeta_\sigma\,\frac{\sqrt{1+p_\sigma^2+q_\sigma^2}}{W} \\[0.6cm]
  &   & \negthickspace\displaystyle
        -\,\sum_{\sigma=1}^{n-2}
           \zeta_\sigma\,
           \mbox{div}
           \left(
             \frac{p_\sigma q^2-q_\sigma(p\cdot q^t)}{W}\,,
             \frac{q_\sigma p^2-p_\sigma(p\cdot q^t)}{W}
           \right) \\[0.6cm]
  &   & \negthickspace\displaystyle
        -\,\sum_{\sigma=1}^{n-2}
           \left\{
             \frac{p_\sigma q^2-q_\sigma(p\cdot q^t)}{\Gamma W}\,\zeta_\sigma\Gamma_x
             +\frac{q_\sigma p^2-p_\sigma(p\cdot q^t)}{\Gamma W}\,\zeta_\sigma\Gamma_y
           \right\}.
\end{array}
\end{equation}
Note that in the second line
\begin{equation}\label{4.26}
\begin{array}{l}
  \displaystyle
  \mbox{div}
  \left(
    \frac{p_\sigma q^2-q_\sigma(p\cdot q^t)}{W}\,\zeta_\sigma,
    \frac{q_\sigma p^2-p_\sigma(p\cdot q^t)}{W}\,\zeta_\sigma
  \right) \\[0.6cm]
  \hspace*{0.6cm}\displaystyle
  =\,\zeta_\sigma\,
     \mbox{div}
     \left(
       \frac{p_\sigma q^2-q_\sigma(p\cdot q^t)}{W}\,,
       \frac{q_\sigma p^2-p_\sigma(p\cdot q^t)}{W}
     \right) \\[0.6cm]
  \hspace*{1.2cm}\displaystyle
     +\,\frac{p_\sigma^2q^2-p_\sigma q_\sigma(p\cdot q^t)}{W}
     +\frac{q_\sigma^2p^2-p_\sigma q_\sigma(p\cdot q^t)}{W}\,,
\end{array}
\end{equation}
and adding up brings
\begin{equation}\label{4.27}
\begin{array}{l}
  \displaystyle
  \sum_{\sigma=1}^{n-2}
  \zeta_\sigma\,
  \mbox{div}
  \left(
    \frac{p_\sigma q^2-q_\sigma(p\cdot q^t)}{W}\,,
    \frac{q_\sigma p^2-p_\sigma(p\cdot q^t)}{W}
  \right) \\[0.6cm]
  \hspace*{0.3cm}\displaystyle
  =\,\sum_{\sigma=1}^{n-2}
     \mbox{div}
     \left(
       \frac{p_\sigma q^2-q_\sigma(p\cdot q^t)}{W}\,\zeta_\sigma,
       \frac{q_\sigma p^2-p_\sigma(p\cdot q^t)}{W}\,\zeta_\sigma
     \right)
     -\frac{2}{W}\,\Big\{p^2q^2-(p\cdot q^t)^2\Big\}.
\end{array}
\end{equation}
Now, (\ref{4.24}) can be written as
\begin{equation}\label{4.28}
\begin{array}{lll}
  W\negthickspace
  & = & \negthickspace\displaystyle
        \frac{1}{W}
        +\sum_{\sigma=1}^{m-2}\mbox{div}\,\frac{\zeta_\sigma\nabla\zeta_\sigma}{W}
        -\frac{p^2q^2-(p\cdot q^t)^2}{W}
        -2\sum_{\sigma=1}^{n-2}
          H(X,\widetilde N_\sigma)\zeta_\sigma\,\frac{\sqrt{1+p_\sigma^2+q_\sigma^2}}{W} \\[0.6cm]
  &   & \negthickspace\displaystyle
        +\,\sum_{\sigma=1}^{n-2}
           \mbox{div}
           \left(
             \frac{p_\sigma q^2-q_\sigma(p\cdot q^t)}{W}\,\zeta_\sigma,
             \frac{q_\sigma p^2-p_\sigma(p\cdot q^t)}{W}\,\zeta_\sigma
           \right) \\[0.6cm]
  &   & \negthickspace\displaystyle
        +\,\sum_{\sigma=1}^{n-2}
           \left\{
             \frac{p_\sigma q^2-q_\sigma(p\cdot q^t)}{\Gamma W}\,\zeta_\sigma\Gamma_x
             +\frac{q_\sigma p^2-p_\sigma(p\cdot q^t)}{\Gamma W}\,\zeta_\sigma\Gamma_y
           \right\}.
\end{array}
\end{equation}
\item[3.]
This last identity will be integrated by parts.
\begin{itemize}
\item[(i)]
First, observe that due to $\frac{1}{W}\,\le 1,$ if holds $\int\hspace{-0.15cm}\int_{\Omega}\frac{1}{W}\,dxdy\le{\mathcal A}\,[\Omega].$
\item[(ii)]
The second is evaluated as
\begin{equation}\label{4.29}
  \sum_{\sigma=1}^{n-2}\,
  \int\hspace{-0.25cm}\int\limits_{\hspace{-0.3cm}\Omega}
  \mbox{div}\,\frac{\zeta_\sigma\nabla\zeta_\sigma}{W}\,dxdy
  \le\sum_{\sigma=1}^{n-2}\,
     \int\limits_{\partial\Omega}
     \frac{|\nabla\zeta_\sigma\cdot\nu^t|}{W}\,|\zeta_\sigma|\,ds
  \le(n-2){\mathcal L}[\partial\Omega]\,\max_{\sigma}\,\|\zeta_\sigma\|_{0,\partial\Omega}
\end{equation}
taking $\frac{\sqrt{1+p_\sigma^2+q_\sigma^2}}{W}\le 1$ for $\sigma=1,\ldots,n-2$ into account.
\item[(iii)]
The third term is non-positive by H\"older's inequality.
\item[(iv)]
Analogously, we have
\begin{equation}\label{4.30}
  \hspace*{-0.3cm}
  2\sum_{\sigma=1}^{n-2}
   \int\hspace{-0.25cm}\int\limits_{\hspace{-0.3cm}\Omega}
   H(X,\widetilde N_\sigma)\zeta_\sigma\,\frac{\sqrt{1+p_\sigma^2+q_\sigma^2}}{W}\,dxdy
  \le 2(n-2)h_0\,{\mathcal A}[\Omega]\,\max_{\sigma}\|\zeta_\sigma\|_{0,\Omega}\,.
\end{equation}
\item[(v)]
We consider the second line in (\ref{4.28}): First, note that
\begin{equation}\label{4.31}
  p_\sigma q^2-q_\sigma(p\cdot q^t)
  =\sum_{\vartheta=1}^{n-2}(p_\sigma q_\vartheta-p_\vartheta q_\sigma)q_\vartheta\,,\quad
  q_\sigma p^2-p_\sigma(p\cdot q^t)
  =\sum_{\vartheta=1}^{n-2}(q_\sigma p_\vartheta-q_\vartheta p_\sigma)p_\vartheta\,,
\end{equation}
and, therefore, multiplication by $p_\sigma$ resp. $q_\sigma,$ and summation brings
\begin{equation}\label{4.32}
  p^2q^2-(p\cdot q^t)^2
  =\frac{1}{2}\sum_{\sigma,\theta=1}^{n-2}(p_\sigma q_\theta-p_\theta q_\sigma)^2\,.
\end{equation}
Integration yields
\begin{equation}\label{4.33}
\begin{array}{l}
  \displaystyle
  \sum_{\sigma=1}^{n-2}\,
  \int\hspace{-0.25cm}\int\limits_{\hspace{-0.3cm}\Omega}
  \mbox{div}
  \left(
    \frac{p_\sigma q^2-q_\sigma(p\cdot q^t)}{W}\,\zeta_\sigma,
    \frac{q_\sigma p^2-p_\sigma(p\cdot q^t)}{W}\,\zeta_\sigma
  \right)dxdy \\[0.8cm]
  \hspace*{0.6cm}\displaystyle
  =\,\sum_{\sigma,\omega=1}^{n-2}\,
     \int\hspace{-0.25cm}\int\limits_{\hspace{-0.3cm}\Omega}
     \mbox{div}
     \left(
       \frac{(p_\sigma q_\omega-p_\omega q_\sigma)q_\omega}{W}\,\zeta_\sigma,
       -\,\frac{(p_\sigma q_\omega-p_\omega q_\sigma)p_\omega}{W}\,\zeta_\sigma
     \right)dxdy \\[0.8cm]
  \hspace*{0.6cm}\displaystyle
  =\,\sum_{\sigma,\omega=1}^{n-2}\,
     \int\limits_{\partial\Omega}
     \frac{(p_\sigma q_\omega-p_\omega q_\sigma)\zeta_\sigma}{W}\,
     (q_\omega,-p_\omega)\cdot\nu^t\,ds \\[0.8cm]
  \hspace*{0.6cm}\displaystyle
  \le\,\sum_{\sigma,\omega=1}^{n-2}
       \int\limits_{\partial\Omega}
       |\zeta_\sigma||D^\top\zeta_\omega|\,ds \\[0.8cm]
  \hspace*{0.6cm}\displaystyle
  \le\,(n-2)^2\,{\mathcal L}[\partial\Omega]\,
       \max_{\sigma}
       \|\zeta_\sigma\|_{0,\partial\Omega}\||D^\top\zeta_\sigma\|_{0,\partial\Omega}\,.
\end{array}
\end{equation}
\item[(vi)]
To control the last term in (\ref{4.28}) we use again (\ref{4.30}), (\ref{4.31}), and $(p_\sigma q_\vartheta-p_\vartheta q_\sigma)^2\le 2[p^2q^2-(p\cdot q^t)^2]$ from (\ref{4.32})
\begin{equation}\label{4.34}
  \frac{|p_\sigma q^2-q_\sigma(p\cdot q^t)|}{W}
  \le\sum_{\vartheta=1}^{n-2}
     \frac{|p_\sigma q_\vartheta-q_\sigma p_\vartheta||q_\vartheta|}{W}
     \,\le\,\sqrt{2}\,
     \sum_{\vartheta=1}^{n-2}\,
     \frac{\sqrt{p^2q^2-(p\cdot q^t)^2}}{W}\,|q_\vartheta|,
\end{equation}
and analogously for $\frac{|q_\sigma p^2-p_\sigma(p\cdot q^t)}{W},$ and we obtain
\begin{equation}\label{4.35}
  \frac{|p_\sigma q^2-q_\sigma(p\cdot q^t)|}{\Gamma W}\le\sqrt{2}\,(n-2)\Gamma_0^{-1}|q|,\quad
  \frac{|q_\sigma p^2-p_\sigma(p\cdot q^t)|}{\Gamma W}\le\sqrt{2}\,(n-2)\Gamma_0^{-1}|p|.
\end{equation}
Then we may estimate
\begin{equation}\label{4.36}
\hspace*{-1cm}
\begin{array}{l}
  \displaystyle
  \sum_{\sigma=1}^{n-2}\,
  \int\hspace{-0.25cm}\int\limits_{\hspace{-0.3cm}\Omega}
  \left\{
    \frac{p_\sigma q^2-q_\sigma(p\cdot q^t)}{W}\,\zeta_\sigma\Gamma_x
    +\frac{q_\sigma p^2-p_\sigma(p\cdot q^t)}{W}\,\zeta_\sigma\Gamma_y
  \right\}\,dxdy \\[0.8cm]
  \hspace*{0.6cm}\displaystyle
  \le\,\sum_{\sigma=1}^{n-2}
       \Gamma_2
       \int\hspace{-0.25cm}\int\limits_{\hspace{-0.3cm}\Omega}
       |\zeta_\sigma|
       \left\{
         \frac{|p_\sigma q^2-q_\sigma(p\cdot q^t)|}{W}
         +\frac{|q_\sigma p^2-p_\sigma(p\cdot q^t)|}{W}
       \right\}\,dxdy \\[0.8cm]
  \hspace*{0.6cm}\displaystyle
  \le\,\frac{\sqrt{2}\,(n-2)^2\Gamma_2}{\Gamma_0}\,\max_{\sigma}\|\zeta_\sigma\|_{0,\Omega}
       \int\hspace{-0.25cm}\int\limits_{\hspace{-0.3cm}\Omega}
       \sqrt{1+p^2+q^2+p^2q^2-(p\cdot q^t)^2}\,dxdy.
\end{array}
\end{equation}
\end{itemize}
\item[6.]
Taking our results together, we arrive at
\begin{equation}\label{4.37}
\hspace*{-0.6cm}
\begin{array}{lll}
  {\mathcal A}[\zeta]\negthickspace
  & \le & \negthickspace\displaystyle
          {\mathcal A}[\Omega]
          +(n-2){\mathcal L}[\partial\Omega]\max_{\sigma}\|\zeta_\sigma\|_{0,\partial\Omega}
          +2(n-2)h_0\,{\mathcal A}[\Omega]\max_{\sigma}\|\zeta_\sigma\|_{C^0(\Omega)} \\[0.4cm]
  &     & \negthickspace\displaystyle
          +\,(n-2)^2\,{\mathcal L}[\partial\Omega]
             \max_{\sigma}
             \|\zeta_\sigma\|_{0,\partial\Omega}\|D^\top\zeta_\sigma\|_{0,\partial\Omega} \\[0.4cm]
  &     & \negthickspace\displaystyle
          +\,\frac{\sqrt{2}\,(n-2)^2\Gamma_2}{\Gamma_0}\,
             \max_{\sigma}\|\zeta_\sigma\|_{0,\Omega}\,
             {\mathcal A}[\zeta].
\end{array}
\end{equation}
Rearranging proves the statement.
\end{itemize}
\end{proof}
\begin{remark}
In (\ref{4.36}) we need the special form $\Gamma=\Gamma(x,y).$ Otherwise, due to (\ref{4.35}) there would remain terms quadratically in $p_\vartheta,$ $q_\vartheta$ in the integrand.
\end{remark}
\subsection{Minimal surfaces}
We consider the special case $\Gamma\equiv 1$ (that is, $h_0=0$):
\begin{corollary}
Let $\zeta\in C^2(\Omega,\mathbb R^{n-2})\cap C^1(\overline\Omega,\mathbb R^{n-2})$ solve the minimal surface system (\ref{4.5}). Then
\begin{equation}\label{4.38}
\hspace*{-0.6cm}
\begin{array}{lll}
  {\mathcal A}[\zeta]\negthickspace
  & \le & \negthickspace\displaystyle
          {\mathcal A}[\Omega]
          +(n-2){\mathcal L}[\partial\Omega]\max_{\sigma=1,\ldots,n-2}\|\zeta_\sigma\|_{0,\partial\Omega} \\[0.4cm]
  &     & \negthickspace\displaystyle
          +\,(n-2)^2{\mathcal L}[\partial\Omega]
             \max_{\sigma=1,\ldots,n-2}
             \|\zeta_\sigma\|_{0,\partial\Omega}\|D^\top\zeta_\sigma\|_{0,\partial\Omega}\,.
\end{array}
\end{equation}
\end{corollary}

\vspace*{0.8cm}
\noindent
Steffen Fr\"ohlich\\
Freie Universit\"at Berlin\\
Fachbereich Mathematik und Informatik\\
Institut f\"ur Mathematik I\\
Arnimalle 2-6\\
D-14195 Berlin\\
Germany\\[0.2cm]
e-mail: sfroehli@mi.fu-berlin.de

\end{document}